\documentclass{ifacconf}

\usepackage{graphicx}      
\usepackage{natbib}        
\usepackage{amsmath}
\usepackage{amssymb}
\usepackage{subcaption}
\usepackage[dvipsnames]{xcolor}
\usepackage{algorithm}
\usepackage{algpseudocode}
\usepackage{comment}
\usepackage{xurl}
\algrenewcommand\algorithmicrequire{\textbf{Input:}}

\newcommand\Pdaqp{\mathcal{P}_\text{DAQP}}
\newcommand\Ptinympc{\mathcal{P}_\text{TinyMPC}}

\begin{document}
\begin{frontmatter}

\title{On Real-Time Feasibility of High-Rate MPC using an Active-Set Method on Nano-Quadcopters} 

\thanks[footnoteinfo]{This work was partially supported by the Wallenberg AI, Autonomous
Systems and Software Program (WASP) funded by the Knut and Alice
Wallenberg Foundation}

\author[First]{J. Wikner} 
\author[Second]{D. Arnström} 
\author[Third]{D. Axehill}

\address[First]{Division of Automatic Control, 
   Linköping University, Linköping, Sweden (e-mail: joel.wikner@liu.se)}
\address[Second]{(e-mail: daniel.arnstrom@gmail.com)}
\address[Third]{Division of Automatic Control, 
   Linköping University, Linköping, Sweden, (e-mail: daniel.axehill@liu.se)}

\begin{abstract}                
Deploying Model Predictive Control on nano-scale aerial platforms is challenging due to the severe computational limitations of onboard microcontrollers. This paper presents an experimental study with computational focus of a dual active-set solver (DAQP) applied to the low-level control of a Crazyflie 2.1 nano-quadcopter. Unlike previous approaches that utilize first-order methods for low-level stabilization or restrict active-set solvers to high-level control, this work demonstrates the successful deployment of a dual active-set method executing at 500 Hz on an STM32F405 microcontroller. We provide a direct benchmark against a state-of-the-art ADMM-based solver (TinyMPC), revealing that the active-set method yields lower execution times for the flight envelopes considered. Furthermore, to manage the complexity of certifying the solver's real-time feasibility, we introduce a data-driven set selection method using Principal Component Analysis. Using this approach, state-of-the-art real-time certification methods can be applied to confirm or reject real-time feasibility already off-line prior to flight. More generally, the experimental results validate that active-set methods are a highly competitive alternative for high-frequency control on resource-constrained hardware and the possibility for their real-time certification adds reliability.
\end{abstract}

\begin{keyword}
Applications of optimal control, real-time optimal control, nano-quadcopter
\end{keyword}

\end{frontmatter}

\section{Introduction}
Model predictive control (MPC) has become known as a highly successful control paradigm for complex systems, especially in the process industry. It is conceptually simple, and its ability to handle hard control constraints make it a useful control-methodology in practice \citep{mayne2014model,borrelli2017predictive}. High-performance control is paramount for enabling tiny, agile robots, such as nano-quadcopters, to operate autonomously in complex and constrained environments \citep{mellinger2012trajectory,johansen2014toward}. MPC is the preferred technique also for these applications, as its ability to explicitly handle state and input constraints is essential for safe and aggressive maneuvers. Unstable systems such as drones pose additional challenges, as they require reliable stabilizing control actions to function. However, achieving such performance onboard of tiny platforms can be challenging due to the limited computational resources onboard, in particular if models including many states are used \citep{mellinger2012trajectory}.

Due to this, model-predictive control (MPC) is typically not preferred when using online onboard control in such scenarios, as the computational cost for MPC can be high \citep{nguyen2021model}. Nevertheless, recent research \citep{nguyen2024tinympc,linder2024real} has shown that there is some hope for utilizing online real-time MPC for resource constrained microcontroller units (MCUs). This has been achieved by in practice using state-of-the-art quadratic programming (QP) solvers to solve the MPC problem. In particular, it has been shown that Linear MPC can be used to control the nano-quadcopter Crazyflie developed by Bitcraze at different levels. Still, applications of efficient second-order active-set methods have not yet been reported for the challenging inner loop of controlling a 12-state nano-quadcopter system. Furthermore, real-time certification has never been reported for MPC applied to this platform.

A strongly related work is the one recently reported in \cite{nguyen2024tinympc}, where an alternating direction method of multipliers (ADMM) algorithm is presented for the nano-quadcopter control application. The algorithm is accelerated using MPC problem structure to precompute and cache matrices to avoid expensive matrix inversions when solving convex QPs online. Exploiting MPC problem structure yields a low memory footprint, allowing the TinyMPC solver to be deployed on an STM32F405 Cortex-M4 MCU clocked at 168MHz to perform low-level control of the Crazyflie 2.1.

During the same period, \cite{linder2024real} independently deployed the dual active-set QP solver DAQP on the Crazyflie to perform high-level control, using the MPC-controller online to compute setpoints for the Bitcraze developed PID-controller.

A second focus in this work is to make use of recent advances in complexity certification \citep{cimini2017exact,arnstrom2021unifying} and worst-case execution time (WCET) analysis \citep{arnstrom2023exact} for active-set quadratic programming algorithms and software, and investigate if and how such tools can be used in this time-critical application. The thesis of \cite{linder2024real} showed the real-time capabilities of the DAQP solver executed on the STM32F405 using the tools developed in \cite{arnstrom2023exact}, concluding that the 100Hz control update frequency was sufficient to achieve the required control performance in the experiments considered. However, no results for the challenging inner loop has so far been reported.

The contribution in this paper are:
\begin{itemize}
    \item Implementation of a \textit{low-level} MPC-controller built on a model containing $12$ states, updated at $500\,$Hz for the Crazyflie nano quadcopter. The MPC-problem is solved onboard online using the dual active-set solver DAQP, and this work represents, to the authors' best knowledge, the first demonstration of that it is indeed possible to use a dual active-set solver on a resource constrained microcontroller for low-level control of a nano quadcopter.
   \item A qualitative experimental comparison of the practical computational performance of the ADMM-based solver TinyMPC and the second-order dual active-set solver DAQP on the STM32F405 Cortex-M4 microcontroller on the Crazyflie is performed. The result illustrates that the more complicated active-set solver is consistently a faster choice more often able to meet the real-time requirements also in this resource-constrained application, as compared to TinyMPC.
   \item The WCET analysis framework recently introduced in \cite{arnstrom2023exact} is applied to analyze the real-time properties of the DAQP solver. This gives novel insight how this tool can be used in a challenging at-the-resource-limit control application to accept and reject designs already at the design stage. The tool is compared with traditional sampling of the state space, which it is shown to outperform.
   \item A novel data-driven strategy based on principal component analysis (PCA) is introduced to increase the off-line computational efficiency of the real-time certification procedure from \cite{arnstrom2023exact}. It also narrows down the size of the parameter set for certification, which turns out to be relevant in at-the-resource-limit control applications. Using this procedure, DAQP turns out not only to be the fastest solver in the comparison, it can also be a priori real-time certified.
\end{itemize}

We would like to stress that while the control performance is implicitly evaluated in our experiments, the main focus in this work is to investigate the real-time computational properties to set the stage for future research more focused on control performance.

The paper is organized as follows: Section 2 gives an overview of the problem formulation, Section 3 describes the main contributions of the paper, Section 4 provides a summary of the experiments and Section 5 presents the conclusions.

\section{Problem Description}

\subsection{Quadcopter modeling}
A nano-qaudcopter is a rigid-body aerial vehicle, actuated by four rotors arranged in a cross configuration. Control is achieved by modulating rotor thrusts and torques individually, to generate desired force and torque. A body-fixed frame is used to describe thrust directions and rotational motion given by attitude $\eta = [\phi,\theta,\psi]^T$ with angular velocity $\omega = [p,q,r]^T$, and an inertial frame is used to describe position $p = [x,y,z]^T$ and velocity $V = [\dot x, \dot y, \dot z]^T$ \cite{mahony2012multirotor}. The state vector is given by
\begin{equation}
    z = [p,V,\eta,\omega]^T.
\end{equation}

\subsection{Linear model predictive control}
Assuming linear time invariant (LTI) dynamics defined by matrices $F\in\mathbb{R}^{n_z \times n_z}$ and $G\in\mathbb{R}^{n_z \times n_u}$, quadratic cost with stage cost penalty matrices $Q\in \mathbb{S}^{n_z}_{+}, R \in \mathbb{S}^{n_u}_{++}$ and terminal cost $P\in\mathbb{S}^{n_z}_{+}$, as well as polyhedral constraints with matrices $A_z\in\mathbb{R}^{n_c\times n_z}, A_u\in \mathbb{R}^{n_c\times n_u}$ and $b\in \mathbb{R}^{n_c}$, the optimal control problem (OCP) 
\begin{equation}\label{eq:mpc}
\begin{aligned}
\min_{u_0,\dots,u_{N-1}} \quad &  \frac{1}{2}z_N^\top P z_N + \frac{1}{2}\sum_{k=0}^{N-1} \left( z_k^\top Q z_k + u_k^\top R u_k \right) \\
\text{subject to} \quad & z_0 = \hat{z}_0\\
& z_{k+1} = F z_k + G u_k, \quad k=0,\dots,N-1,\\
& A_zz_k + A_u u_k \le b_u, \quad k=0,\dots,N-1,\\
& A_f z_N \le b_f
\end{aligned}
\end{equation}
is in linear MPC solved at each time-step in a receding horizon fashion. The first computed control signal $u_0$ is the one actually applied as input to the system. Here $z_k\in\mathbb{R}^{n_z}, u_k\in\mathbb{R}^{n_u}$ are states and inputs respectively, and $A_f\in \mathbb{R}^{n_f\times n_z}$, $b\in \mathbb{R}^{n_f}$ define the terminal set. The initial state $z_0$ in the finite-time horizon is set to the current measured or estimated system state $\hat{z}_0$, closing the loop.

\subsection{Quadratic programming}
Given $\hat{z}_0$, the OCP in \eqref{eq:mpc} can be cast as a QP using a condensed formulation, see for instance \cite{borrelli2017predictive}. The primal and dual formulations of a general QP are given as

\begin{center}
\fbox{%
  \begin{minipage}[t]{0.405\linewidth}
    \centering \textbf{Primal QP}\\[6pt]
    \[
    \begin{aligned}
      \min_{x}\hspace{0.5em} & \tfrac{1}{2} x^\top H x + f^\top x\\
      \text{s.t.}\hspace{0.5em}& [A]_i x \le [b]_i, \forall i\in\mathbb{N}_m
    \end{aligned}
    \]
  \end{minipage}
}
\hspace{0.02\linewidth}
\fbox{%
  \begin{minipage}[t]{0.405\linewidth}
    \centering \textbf{Dual QP}\\[6pt]
    \[
    \begin{aligned}
      \min_{\lambda}\hspace{0.5em} & \tfrac{1}{2}\lambda^\top MM^T\lambda + d^\top \lambda\\
      \text{s.t.}\hspace{0.5em} & [\lambda]_i \ge 0,\; \forall i\in\mathbb{N}_m
    \end{aligned}
    \]
  \end{minipage}
}
\end{center}

with unique primal and dual solutions $x^{\star}, \lambda^{\star}$, respectively, for convex QPs with $H\in \mathbb{S}_{++}^{n}$ and $f\in \mathbb{R}^n$. Moreover, the feasible set is given by $A\in\mathbb{R}^{m\times n}$ and $b\in\mathbb{R}^m$ where $[.]_i$ denotes the \textit{i}th row of a matrix. The dual problem is defined by $MM^T=(A H^{-1} A^\top)$ and $d=(b + A H^{-1} f)$. \cite{arnstrom2023real}

Common methods for finding $x^{\star}$ are active-set methods \citep{goldfarb1983numerically,schmid1994quadratic,axehill2006mixed,axehill2008dual,ferreau2014qpoases,saraf2019bounded,arnstrom2022dual}, operator splitting methods \citep{stellato2020osqp} and interior-point methods \cite{wang2009fast,frison2020hpipm}. 

Much of the complexity in solving QPs comes from the inequality constraints $[A]_ix \leq [b]_i$. Active-set methods approach this problem by solving a sequence of equality constrained QPs (EQP) in the form

\begin{equation}\label{eq:eqp}
\begin{aligned}
\min_{x} \quad & \frac{1}{2}x^\top H x +f^{T}x\\
\text{subject to} \quad & [A]_ix = [b]_i, \quad \forall i\in \mathcal{A}^{\star}\\
\end{aligned}
\end{equation}
until a solution to \eqref{eq:eqp} satisfies the Karush Kuhn Tucker (\textit{KKT}) necessary and sufficient optimality conditions 
\begin{align}
H x^\star + A^\top \lambda^\star &= -f \\
Ax^\star &\le b \\
\lambda^\star &\ge 0 \\
\lambda_i^\star ([A]_i x^\star - [b]_i) &= 0, \quad i = 1, \dots, m.&
\end{align}
of the original problem. Hence, a large part of the active-set algorithm reduces to a combinatorial problem of finding the optimal active-set $\mathcal{A}^{\star}= \{\,i\in \mathbb{N}_m \, \vert \, [A]_ix_{\star}=[b]_i\, \}$, \citep{arnstrom2023real}.

Essentially, the active-set algorithm starts with a working set $\mathcal{W} \subseteq \mathbb{N}_m$ initialized to $\mathcal{W}_0$, iteratively adds candidates from $\mathbb{N}_m\setminus\mathcal{W}$ or removes elements from $\mathcal{W}$, solves the KKT system
\begin{equation}\label{eq:KKT_system}
\begin{pmatrix}
H & [A]_\mathcal{W}^\top \\
[A]_\mathcal{W} & 0
\end{pmatrix}
\begin{pmatrix}
x \\ [\lambda]_\mathcal{W}
\end{pmatrix}
=
\begin{pmatrix}
- f \\ [b]_\mathcal{W}
\end{pmatrix}
\end{equation}
and checks for optimality. An overview of the active-set algorithm is given in Algorithm \ref{alg:active-set}, and a complete description can be found in \cite{wright1999numerical}.

\begin{algorithm}
\caption{Archetypal active-set QP method}
\label{alg:active-set}
\begin{algorithmic}[1]
\Require $H$, $f$, $A$, $b$, $\mathcal{W}_0$, $\lambda_0$
\State $\mathcal{W} \leftarrow \mathcal{W}_0$
\State $\lambda \leftarrow \lambda_0$
\Repeat
    \State $(x,\lambda) \leftarrow \text{Solve KKT system in \eqref{eq:KKT_system} given } \mathcal{W}$
    \If{$(x,\lambda)$ is primal and dual feasible}
        \State \Return $(x^\star, \lambda^\star, \mathcal{A}^\star) \leftarrow (x,\lambda,\mathcal{W})$
    \Else
        \State Add $i\in \mathbb{N}_m \setminus \mathcal{W} $ to $\mathcal{W}$ or remove element from  $\mathcal{W}$ based on primal and/or dual violation of $x$ and $\lambda$
    \EndIf
\Until{termination}
\end{algorithmic}
\end{algorithm}

As Algorithm \ref{alg:active-set} computes $x^{\star}$ iteratively, it also generates a working-set sequence $\{ \mathcal{W}_j \}_{j=0}^{k}$ where $k$ denotes the number of iterations at which the algorithm terminates. Since the dual formulation of the QP has the same principle form as the primal, Algorithm \ref{alg:active-set} can also be applied to the dual QP problem as long as subproblems in the form \ref{eq:eqp} with $H \in\mathbb{S}^{n}_{+}$ can be solved. Consequently, the following theory is also valid for dual active-set solvers such as DAQP, \cite{arnstrom2022dual}.

\subsection{Time-complexity certification for MPC}\label{sec:time_cert}

Another way to look at the linear MPC problem is to, instead of assuming $\hat{z}_0$ to be fixed, allow $\hat{z}_0$ to vary within a given set $\Theta = \{\theta:A_{\theta}\theta \leq b_{\theta}\}$ denoted as the parameter-space, with parameter $\theta = \hat{z}_0$. It is shown in \cite{cimini2017exact,arnstrom2021unifying} that, since the primal QP vary with the parameter $\theta$ via the affine functions $f(\theta) = \overline{f}+f_{\theta}\theta$ and $b(\theta) = \overline{b} + W_{\theta}\theta$, it is possible to \textit{exactly} determine the working set sequence $\{ \mathcal{W(\theta)}_j \}_{j=0}^{\bar{k}(\theta)}$ and the exact number of iterations $\bar{k}(\theta)\in\mathbb{N}$ as a piecewise constant function for $\theta \in \Theta$ using, for example, the method in \cite{arnstrom2021unifying}. 

More specifically, the working set sequence $\mathcal{W}$ and $\bar{k}$ remain fixed within the parameter regions $\{ \Theta_{j} \}_{j=0}^{M}$ in the parameter space that partition $\Theta$, such that $\Theta = \bigcup_{i=1}^M \ \Theta_{j}$ and $\overset{\circ}{\Theta}_i \cap \overset{\circ}{\Theta}_j = \varnothing$ for $i\neq j$ \cite{arnstrom2023exact}. Here, $\overset{\circ}{\Theta}$ denotes the interior of $\Theta$. 

Beyond algorithms, we define a \textit{program} $\mathcal{P}$ as a software implementation of an algorithm that solves \eqref{eq:mpc} on hardware (e.g. an MCU). The behavior of a program depends not only on the algorithm itself, but also on the compiler that generates the machine code and on the hardware that executes the code.

Since the number of regions $\Theta_j$ is finite, the number of working set sequences to analyze is also finite. Hence, it suffices to analyze how the active-set solver behaves for a chosen $\theta \in \Theta_j$ in each region. Consequently, \cite{arnstrom2023exact} show that it is possible to give a WCET certificate for a program $\mathcal{P}$ that implements Algorithm \ref{alg:active-set} on specific hardware, given some implicit requirements on $\mathcal{P}$. See Assumption B.1-4 and Theorem 1 in \cite{arnstrom2023exact}. 

The main idea introduced in \cite{arnstrom2023exact} is to, given a certain target hardware, measure the execution time $\tau_j$ of the program $\mathcal{P}$ used to solve the QPs of interest in carefully selected samples $\theta_j \in \Theta_j$ for $j=1,..,M$, generated by the framework in \cite{arnstrom2021unifying} and compute the (tight) WCET according to $\overline{\tau} = \max\limits_{j} \tau_j$.

Since the Crazyflie uses a real-time OS, we have in this explorative and practically grounded work relaxed some formal assumptions B.1-4 in \cite{arnstrom2023exact} related to, e.g., the processor state. Furthermore, the reduced floating point precision of the MCU has not been explicitly considered here. However, the used framework has in principle support for taking limited precision into account as presented in \cite{arnstrom2022lift}, but such analysis will increase the off-line computational burden and has been considered future work. Still, as will be shown in the experimental section, the analysis will align very well with the observations made in practice.

We investigate the computational complexity of a given program $\mathcal{P}$ using the time complexity certification framework in \cite{arnstrom2021unifying} implemented in the Julia package ASCertain, which specifically analyzes DAQP. The procedure is summarized in Algorithm \ref{alg:cert_data}, which produces execution-time measurements stored in the vector $\tau\in\mathbb{R}^M$ that enable sample-wise comparison between different programs using, for instance, subtraction.

\begin{algorithm}
\caption{Execution time measurements}
\label{alg:cert_data}
\begin{algorithmic}[1]

\Require Parameter set $\Theta$, OCP in \eqref{eq:mpc} and Program $\mathcal{P}$

\State $\{H,A,f(\theta),b(\theta)\}$ $\gets$ Use condensed formulation on OCP \eqref{eq:mpc}
\State $\{ \Theta_{j} \}_{j=0}^{M}$ $\gets$ Use ASCertain on $\{H,A,f(\theta),b(\theta)\}$ for $\Theta$
\State $\mathcal{J} \gets \{1, \ldots, M\}$
\State $\boldsymbol{\tau} \gets [\,]$
\For{$j \in \mathcal{J}$}
    \State $\theta_j \gets$ Select any $\theta_j \in \Theta_j$
    \State $\tau_j \gets$ Measure time for executing $\mathcal{P}$ with $\theta = \theta_j$
    \State $\tau \gets \text{Push}(\tau, \tau_j)$
\EndFor

\State \Return Vector with measured time $\tau \in \mathbb{R}^M$

\end{algorithmic}
\end{algorithm}

Beyond complexity and time certification, the partition from \cite{arnstrom2021unifying} gives parameters $\theta \in \Theta_j$ corresponding to QPs of high relevance when comparing QP solvers. By construction, such an analysis considers all relevant cases, and in particular the hard cases, of interest for the certified solver. Hence, it provides an exact (or more generally upper) bound of the computational complexity for that solver. If an alternative solver is benchmarked in exactly these samples, the observations from the alternative solver will form a lower bound on its worst-case computation time. Hence, such an analysis will be conservative for the certified active-set solver. 

\section{Real-time MPC control of the nano quadcopter using a dual active set solver}\label{sec:realtime}

For the design of the MPC controller given in \eqref{eq:mpc}, we use $Q$, $R$, $N$, $F$, $G$, $A_z$ and $A_u$ provided by the TinyMPC firmware. These matrices define the controller used in the figure-eight experiment in \cite{nguyen2024tinympc}. The terminal cost $Q_f$ was set to the infinite time Riccati solution, and no terminal set was used. This problem data together with the Julia package LinearMPC is used to generate software in C code that implements DAQP. As a result, both TinyMPC and DAQP use the same MPC problem setup, including tuning, as in the TinyMPC firmware available on github at \cite{RoboticExplorationLab_tinympc2025}. From the control results, one can argue that this tuning can be improved, but the decision was made to consistently keep this one since computational performance is the main focus here and to avoid the risk of introducing a bias in the quality of the tunings of the solvers at the limit of their computational capacity. The primal tolerance $\epsilon_p$ is set to 1e-4, and the dual tolerance $\epsilon_d$ is set to 1e-4 for DAQP. The code is compiled using the \texttt{-Ofast} flag in GCC to produce the program $\Pdaqp$. Note that $A_z$ and $A_u$ from the firmware only define input constraints. DAQP uses an empty initial working set $\mathcal{W}_0$ and is not warm-started, as the real-time certification tool ASCertain analyzes DAQP for a given initial working set.


For comparison, the TinyMPC firmware source code was compiled with the \texttt{-Ofast} GCC optimization flag to produce the program $\Ptinympc$. No changes were made to the source code for the TinyMPC solver implementation, apart from raising the maximum iteration limit to 100 to ensure convergence within the designated dual and primal tolerances during runtime.

The states to be controlled in \eqref{eq:mpc} are $z_e=z-z_r$ where $z$ are the position, attitudes, linear velocities and angular velocities provided as estimates by the EKF after the attitudes have been transformed from quarternions to Rodriguez parameters as in \cite{jackson2021planning}. Moreover, $z_r$ is the state reference given as a pre-programmed trajectory stored onboard the Crazyflie. The problem has state dimension $n_z=12$ and input dimension $n_u=4$, where the input $u$ computed by $\Pdaqp$ together with $u_0$, that is nominal input around hover, gives the PWM motor commands 
\begin{equation}\label{eq:motor}
    u_m=u_0+u    
\end{equation}
normalized to the interval $[0,1]$. The MPC controller is executed at 500 Hz with a horizon length of $N=15$.


To evaluate execution time of the QP-solvers with Algorithm \ref{alg:cert_data}, the parameter region $\Theta_b $ was selected as $\Theta_b = \bigl\{\, z \in \mathbb{R}^{12} \;\big|\; 
[\underline{z}]_i +[\underline\epsilon]_i \le [z]_i \le [\overline{z}]_i + [\overline\epsilon]_i, i=1,...,12 \,\bigr\}$, where $\underline{z}_i = \min\limits_{k} [z_{e}^{(k)}]_i$ and $\overline{z}_i = \max\limits_{k} [z_{e}^{(k)}]_i$ were chosen by collecting the measured error states $z_e^{(1)},z_e^{(2)}...,z_e^{(K)}\in\mathbb{R}^{12}$ at 20 Hz during the execution of a test flight in which the drone maintained constant altitude and followed a figure-eight trajectory in the horizontal plane. $\underline\epsilon,\overline\epsilon\in \mathbb{R}^{12}$ were here user defined to, with some margin, cover all operating points if interest. 
\subsection{Data-driven set selection for complexity certification}
Only a limited subset of the states contained in the box $\Theta_b$ considered above is actually reached during flight. An important reason for this is that it is likely that a simple axis-aligned box in the state-space contains many non-physical combinations of states which are practically irrelevant. Hence, often a smaller subset of the initially considered box is sufficient to capture the state trajectories that occur in practice. 

To automatically get a more representative set of the states that are reached during flight, a data-driven approach is used to obtain a hyper-rectangle better representing $z_{e}^{(K)}$ using PCA. For more details on PCA, see \cite{abdi2010principal}. The measured error state are stacked as $Z =[z_e^{(1)}, z_e^{(2)} \cdots z_e^{(K)}]$ with sample mean $\mu = \frac{1}{K} Z \mathbf{1} \in \mathbb{R}^{12}$ where $\mathbf{1}\in\mathbb{R}^{K}$ is a vector that contains only ones. PCA was implemented using singular value decomposition of the centered error states $U\Sigma V^T = Z-\mu$. With $U$ orthonormal, the error states are expressed in a new basis chosen as $\tilde{z}^{(K)} = U^\top (\tilde{z}_e^{(K)} - \mu)$. The upper and lower bounds of the rotated hyper-rectangle are defined as $[\underline{\tilde{z}}]_i = \min\limits_{k} [z_{e}^{(k)}]_i$ and $[\overline{\tilde{z}}]_i = \max\limits_{k} [\tilde{z}_{e}^{(k)}]_i$, respectively. A user-defined constant $\delta$ is introduced to give the possibility for introducing a relative margin $\delta(\overline{\tilde{z}}-\underline{\tilde{z}})$ providing the option to inflate the considered parameter set and hence providing some margin for so far unseen states, resulting in the set $\underline{\tilde{z}}-\delta(\overline{\tilde{z}}-\underline{\tilde{z}}) \leq \tilde{z} \leq \overline{\tilde{z}}+\delta(\overline{\tilde{z}}-\underline{\tilde{z}})$. Inserting $\tilde{z} = U^T(z_e-\mu)$ into the inequality yields 

\begin{equation}
    U_i^\top z_e \le \overline{\tilde{z}} +\delta(\overline{\tilde{z}}-\underline{\tilde{z}}) + U^\top \mu,
\end{equation}
\begin{equation}
    -U_i^\top z_e \le -\underline{\tilde{z}}  + \delta(\overline{\tilde{z}}-\underline{\tilde{z}})- U_i^\top \mu.
\end{equation}

Stacking the inequalities finally gives the rotated hyper-rectangle 

\begin{equation}
A_p = 
\begin{bmatrix}
U^\top \\[4pt]
-U^\top
\end{bmatrix},
\qquad
b_p =
\begin{bmatrix}
\overline{\tilde{z}}+ \delta(\overline{\tilde{z}}-\underline{\tilde{z}}) + U^\top \mu \\[4pt]
-(\underline{\tilde{z}}- \delta(\overline{\tilde{z}}-\underline{\tilde{z}}) + U^\top \mu)
\end{bmatrix},
\end{equation}

with $\Theta_p = \{\, z_e \in \mathbb{R}^n : A_p z_e \le b_p \,\}$.

\section{Experimental Evaluation}
\subsection{Hardware description}

The platform considered in this work is the Bitcraze Crazyflie 2.1 nano-quadcopter, which weighs 29\,g and is equipped with a low-latency radio interface that enables real-time data transfer during flights. The system is tracked with a Qualisys motion capture-system in a $12 \times 12 \times 8$m indoor arena providing high-precision position measurements at 100 Hz to the Crazyflie using the radio interface. Angular velocities and linear accelerations of the drone are measured using an inertial measurement unit (IMU). The MCU on the Crazyflie fuses the Qualisys measurements with the IMU-data with a factory-supplied extended Kalman filter (EKF) to produce state estimates, which are used for control. 

Both the control-algorithms and the EKF are executed on the drone using an STM32F405 Cortex-M4 168 MHz MCU with 196kB of SRAM and 1MB Flash. The experiments performed in this paper were, to our best efforts, started with fully charged stock 250 mAh batteries. No additional Bitcraze expansion deck was attached to the drone, aside from five Qualisys passive markers.

\subsection{Benchmarking procedure}\label{sec:benchmark}

\subsubsection{Parameter sampling strategies}
Two methods were used to generate data for the comparison of TinyMPC and DAQP with respect to the execution time. In the first method, Algorithm \ref{alg:cert_data} with $\Theta_b$ and the OCP specified in Section \ref{sec:realtime} was used for both programs $\Pdaqp$ and $\Ptinympc$. This allowed us to measure the execution times stored in the vectors $\tau_{\Pdaqp}^c$ and $\tau_{\Ptinympc}^c \in \mathbb{R}^{M_c}$ on the MCU for DAQP and TinyMPC, respectively. The upper and lower bounds of $\Theta_b$ were selected as $\overline{z} + \overline\epsilon=[0.6, 1.0, 0.6, 0.5, 0.4, 0.25, 1.6, 3.0, 0.2 , 7.0, 5.0, 0.45]^T$ and $[\underline{z}] +\underline\epsilon =-(\overline{z} + \overline\epsilon)$ respectively.

For the second method, $\{\theta_j^u\}_{j=0}^{M_c}$ was generated from a uniform distribution over the interval $[\underline{z} +\underline\epsilon$, $\overline{z} + \overline\epsilon]$. Note that the same number of data points is generated as in the first method. The execution time for the uniformly distributed data was measured for both programs $\Pdaqp$ and $\Ptinympc$ using Algorithm \ref{alg:cert_data} and replacing Step 4 with $\theta_j = \theta_j^u$, yielding execution times stored in the vectors $\tau_{\Pdaqp}^u$ and $\tau_{\Ptinympc}^u\in \mathbb{R}^{M_c}$ for DAQP and TinyMPC respectively.

Neither TinyMPC nor DAQP were warm-started with the previous solution in these experiments, since the uniform samples are not necessarily close and suitable for warm-starts as in a closed-loop trajectory.


\subsubsection{Microcontroller benchmarks}
\begin{figure}[htbp]
    \centering
    \includegraphics[width=0.9\linewidth]{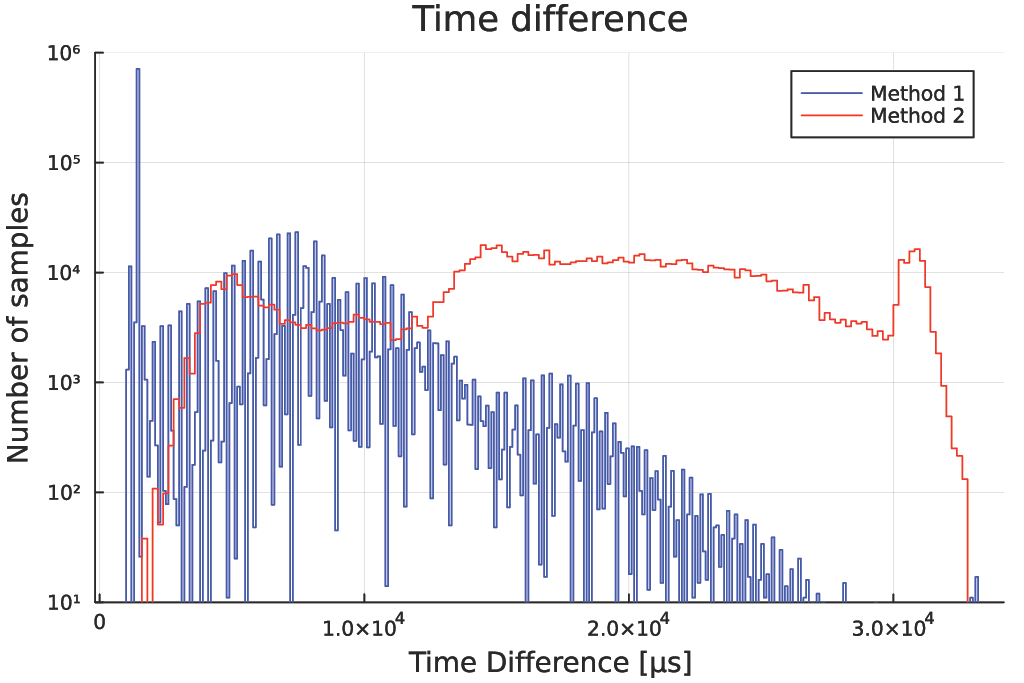}
    \caption{Histogram of difference in execution time distribution $\tau_{\mathcal{P}_{TinyMPC}}^c-\tau_{\mathcal{P}_{DAQP}}^c$ and $\tau_{\mathcal{P}_{TinyMPC}}^u-\tau_{\mathcal{P}_{DAQP}}^u$ for DAQP and TinyMPC. A positive difference implies that DAQP is faster than TinyMPC. All computed differences were positive, demonstrating that DAQP is consistently faster across all tested instances.}
    \label{fig:execution_time}
\end{figure}

The number of data points for both Method 1 and Method 2 is 1179415, and the histogram of the, as elaborated in Section \ref{sec:time_cert}, sample-wise differences $\tau_{\mathcal{P}_{TinyMPC}}^c-\tau_{\mathcal{P}_{DAQP}}^c$ for Method 1 and $\tau_{\mathcal{P}_{TinyMPC}}^u-\tau_{\mathcal{P}_{DAQP}}^u$ for Method 2 is shown in Figure \ref{fig:execution_time}. It can be seen that the difference in computation time is strictly greater than zero, implying that DAQP always computes the solution faster for the parameters considered.

\begin{figure}[htbp]
    \centering
    \includegraphics[width=0.9\linewidth]{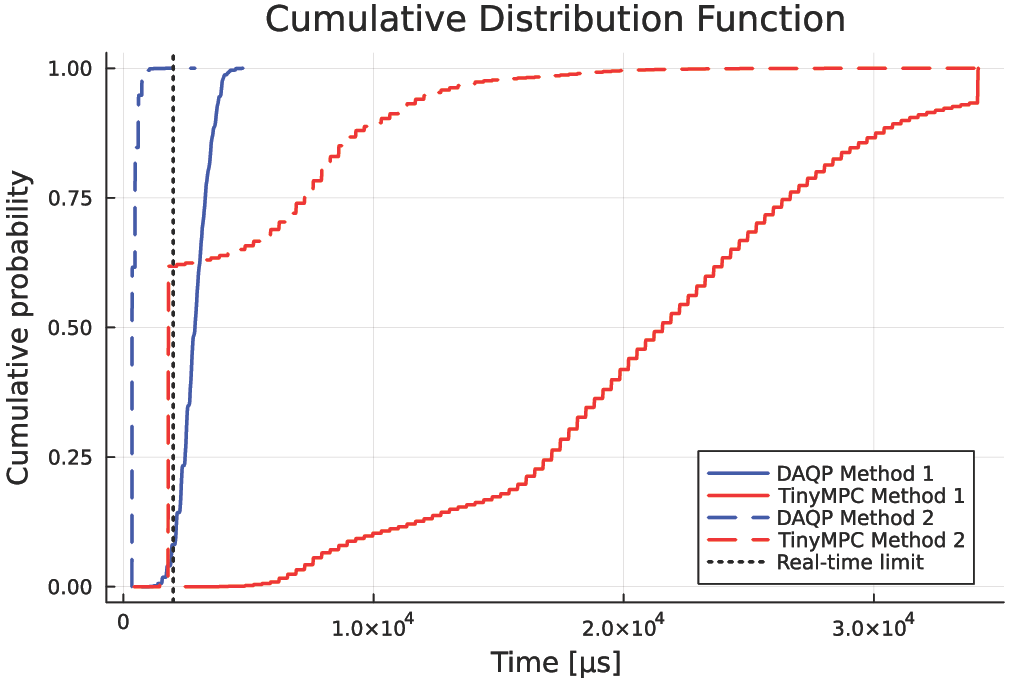}
    \caption{CDFs of execution times $\tau_{\mathcal{P}_{TinyMPC}}^u$ and $\tau_{\mathcal{P}_{DAQP}}^u$ from Method~1 and Method~2. In addition to showing that DAQP outperforms TinyMPC, it shows that Method~1 might give a biased view of the worst-case execution time, which can give a false sense of the real-time capabilities for the controller.}
    \label{fig:PCA}
\end{figure}

Figure \ref{fig:PCA} shows the cumulative distribution function (CDF) of $\tau_{\mathcal{P}_{DAQP}}^u$, $\tau_{\mathcal{P}_{TinyMPC}}^u$, $\tau_{\mathcal{P}_{DAQP}}^c$ and $\tau_{\mathcal{P}_{TinyMPC}}^c$. The plot illustrates that, given the same number of samples, using a uniform distribution to evaluate real-time performance fails to reveal the worst-case scenario for the active-set solver, which may falsely suggest that the solver is suitable for real-time deployment over the flight-envelope in $\Theta_b$.

\subsubsection{Hardware experiments}\label{sec:flight}
To evaluate the real-time capability and control performance of the MPC controllers given in the programs $\mathcal{P}_{DAQP}$ and $\mathcal{P}_{TinyMPC}$, a deliberately minimalistic ("clean") step-response-based trajectory was used. The motivation for the use of this trajectory was that it is fairly simple and can provide gradual challenges in terms of number of activated constraints which is suitable for evaluating the computational performance of the controllers. The trajectory is composed by three steps executed after each other in the $xz$-plane: first one in the positive x-direction, then one in the positive z-direction, and finally one diagonally combining a simultaneous step in the negative x and z directions returning to the start position. Before the drone started the execution of the three maneuvers, it took off and reached hovering state at the starting altitude for the experiments. All non-position states were penalized towards the origin. The result is shown in Figure \ref{fig:flightdaqptinympc}. Here, the $x$-axis defines a certain horizontal direction in the flight arena and the $z$-axis defines the vertical direction pointing upwards. TinyMPC was warm-started with the previous solution in these experiments. While both DAQP and TinyMPC fulfilled the real-time requirement during flight, it can be seen in Figure \ref{fig:flightdaqptinympc} that neither controller manages to stay close to the desired altitude using the tuning from the TinyMPC firmware. This will be discussed more below. The seemingly more complex, but actually less computationally challenging since no constraints are activated, figure-eight trajectory used in~\cite{nguyen2024tinympc} is reproduced and presented as a reference result in Figure~\ref{fig:flightdaqptinympceight}, confirming the similarity of the setups in the two works.

\subsubsection{Simulation study}
\begin{figure}[htbp]
    \centering
    \includegraphics[width=0.9\linewidth]{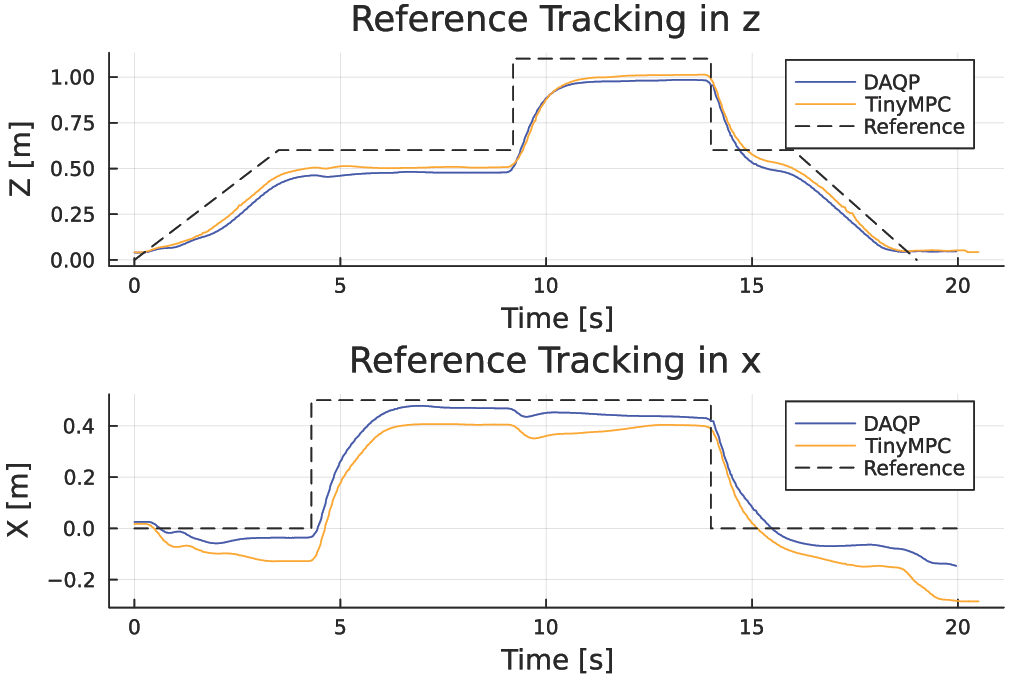}
    \caption{Comparison of waypoint tracking for DAQP vs TinyMPC. The Crazyflie performed an initial takoff, followed by horizontal translation to waypoint 1, ascent to waypoint 2, and return to hover at waypoint 3}
    \label{fig:flightdaqptinympc}
\end{figure}

The MPC controller in program $\mathcal{P}_{DAQP}$ was evaluated in simulation to investigate its performance when applied to simulated nonlinear 6-DOF flight dynamics using code from \cite{RoboticExplorationLab_tinympc2025}. In simulation, the drone executed the same step-based trajectory as in the hardware experiment, and achieved acceptable tracking performance. Interestingly, introducing a small reduction in $u_0$ would yield a static offset in altitude tracking, similar to the tracking offset in Figure \ref{fig:flightdaqptinympc}. This indicates that the altitude-tracking error observed in practice is likely to be caused by that the effect from gravity deviates from the static compensation $u_0$ made according to \eqref{eq:motor}. In a more complete control design, this would have been compensated by integral action. As one would expect, the vertical tracking error can also be reduced by decreasing $R$ in \eqref{eq:mpc}, which then deviates from the tuning suggested in \cite{nguyen2024tinympc}.

\subsubsection{Effect of input weight $R$ on tracking performance} 
\begin{figure}[htbp]
    \centering
    \includegraphics[width=0.9\linewidth]{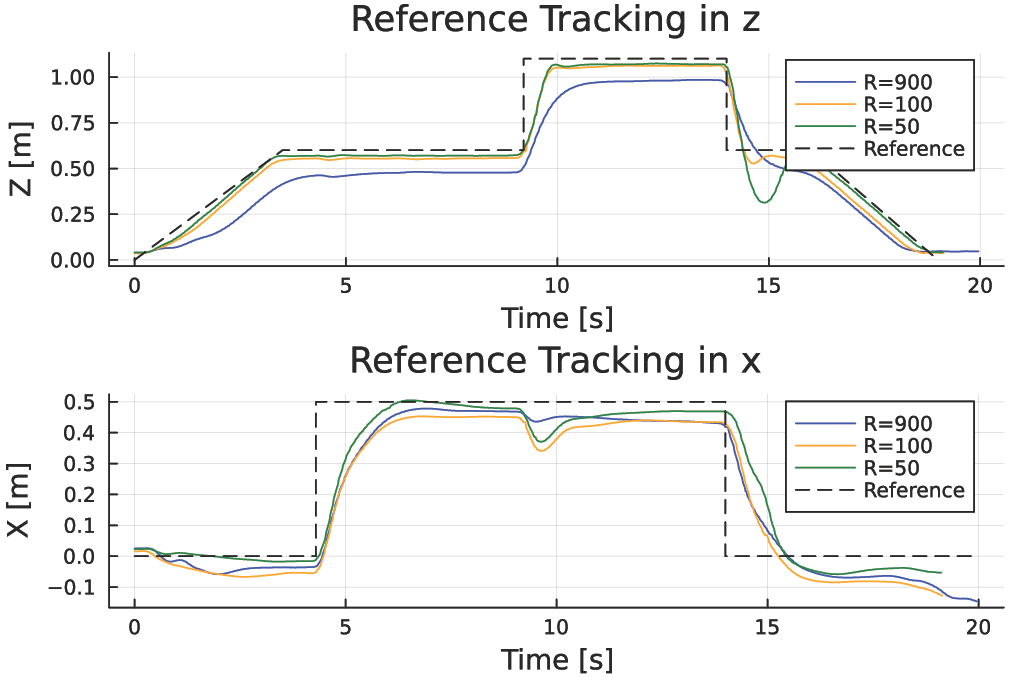}
    \caption{Waypoint tracking performed with DAQP for different tunings on $R$. Decreasing the weight $R$ improves tracking in altitude, but it can eventually damage real-time feasibility.}
    \label{fig:flightdaqp}
\end{figure}

In order to investigate a possible change in behavior, we reduced the weight $R$ from the original tuning of 900 to first 100 and then 50 for all diagonal elements and tested the new tuning on DAQP. The off-diagonal terms remain zero, and the result can be seen in Figure \ref{fig:flightdaqp}. It can be seen as expected that reducing the weight on the control signal generally resulted in a lower steady-state error and improved tracking. However, for the most aggressive choice $50$,  the real-time requirement of $500\,$ Hz was visibly violated during the descent in the last step response diagonally back to the starting position degrading the control performance. This results in the undershoot in the green curve around $15\,$s in Figure~\ref{fig:flightdaqp}.

\subsubsection{Execution time evaluation with PCA-boxes}
Using the measured error states from the flight in Figure \ref{fig:flightdaqp} for the original tuning in the TinyMPC firmware, a tighter hyper-rectangle $\Theta_{pca}^{900}$ is obtained with $\delta$ set to 2. Algorithm \ref{alg:cert_data} was used with $\Theta_{pca}^{900}$, $\mathcal{P}_{DAQP}$ and the OCP specified in Section \ref{sec:realtime}. The worst case execution time was 1375 $\mu $s. 

Hence, using the data-driven approach to obtain $\Theta_{pca}^{900}$ better reflects the observed real-time capability of the controller for online control, rather than using an overly conservative $\Theta_b$, for which DAQP violates the real-time requirement on some states, as seen in Figure \ref{fig:PCA}.

The procedure was repeated for the controller with diagonal elements of $R$ set to 100, to obtain $\Theta_{pca}^{100}$. This time with $\delta$ set to zero. Using Algorithm \ref{alg:cert_data} with the modified controller and $\Theta_{pca}^{100}$ revealed that the controller violates the real-time requirements in the considered parameter set. This demonstrates how Algorithm \eqref{alg:cert_data} can be used to predict when the controller fails to meet hard real-time requirements.

\section{Conclusion}
This work demonstrates that a carefully implemented dual active-set solver can be successfully used on resource-constrained embedded platforms. In particular, high-rate execution of a low-level 12-state MPC controller running entirely onboard on a Crazyflie nano-quadcopter has been successfully demonstrated. Using the WCET analysis framework gives novel insight into real-time feasibility of the controller already at the design stage, and a PCA-based parameter-set selection approach improves the computational tractability and reduces the conservativeness of the certification process. Data for PCA-based set selection can be obtained from both real experiments and from a simulator, highlighting the practicality and versatility of the method.

The results show that dual active-set solvers can deliver both real-time guarantees and high computational performance for nano-quadcopter control, opening up new opportunities for reliable embedded optimization. An experimental comparison against the ADMM-based solver TinyMPC highlights the competitiveness of DAQP, as it consistently is a faster solver, which is of importance for this at-the-resource-limit control application.

Future work includes more focus on the actual control performance (for example, using tools from offset-free MPC) and further improvement of the computational performance of the certification algorithm and implementation.


\begin{ack}
Arvid Linder is acknowledged for the work in his Master's thesis \cite{linder2024real}, providing a good foundation for the work in this paper.
\end{ack}

\section*{DECLARATION OF GENERATIVE AI AND AI-ASSISTED TECHNOLOGIES IN THE WRITING PROCESS}
During the preparation of this work the author(s) used ChatGPT and Gemini in order to assist with text editing and clarity. After using this tool/service, the author(s) reviewed and edited the content as needed and take(s) full responsibility for the content of the publication

\bibliography{ifacconf}             









\appendix
\section{Additional Figures}    
\begin{figure}[htbp]
    \centering
    \includegraphics[width=0.9\linewidth]{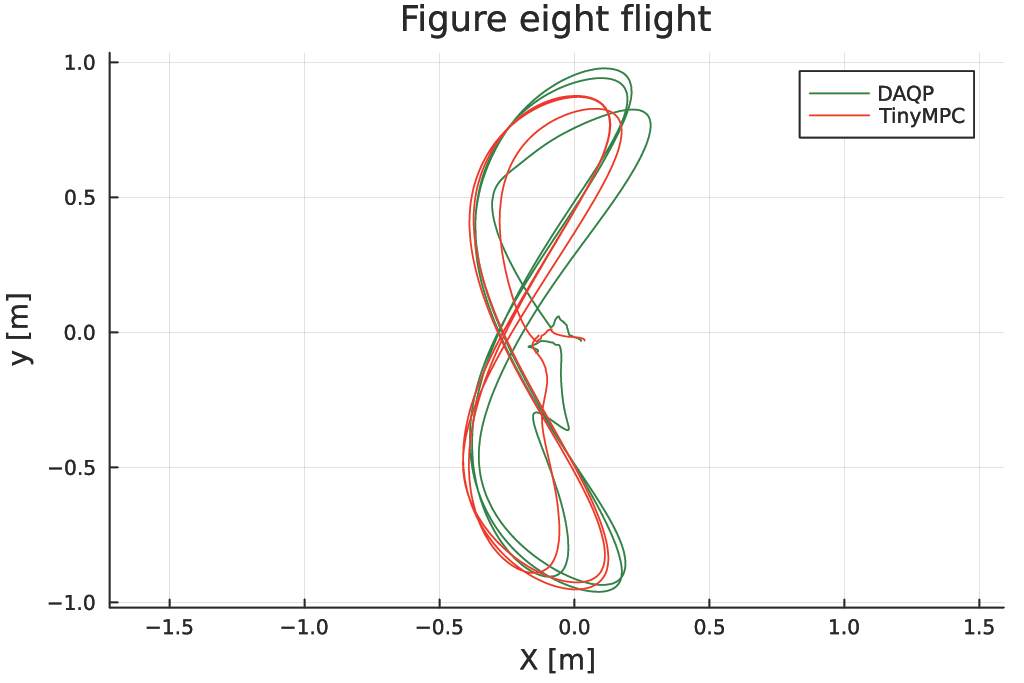}
    \caption{Figure-eight trajectory tracking for DAQP vs TinyMPC. It is noted that this result is well-aligned with what is reported in \cite{nguyen2024tinympc}, confirming that the setup in this work is correct.}
    \label{fig:flightdaqptinympceight}
\end{figure}             
\begin{figure}[htbp]
    \centering
    \includegraphics[width=0.9\linewidth]{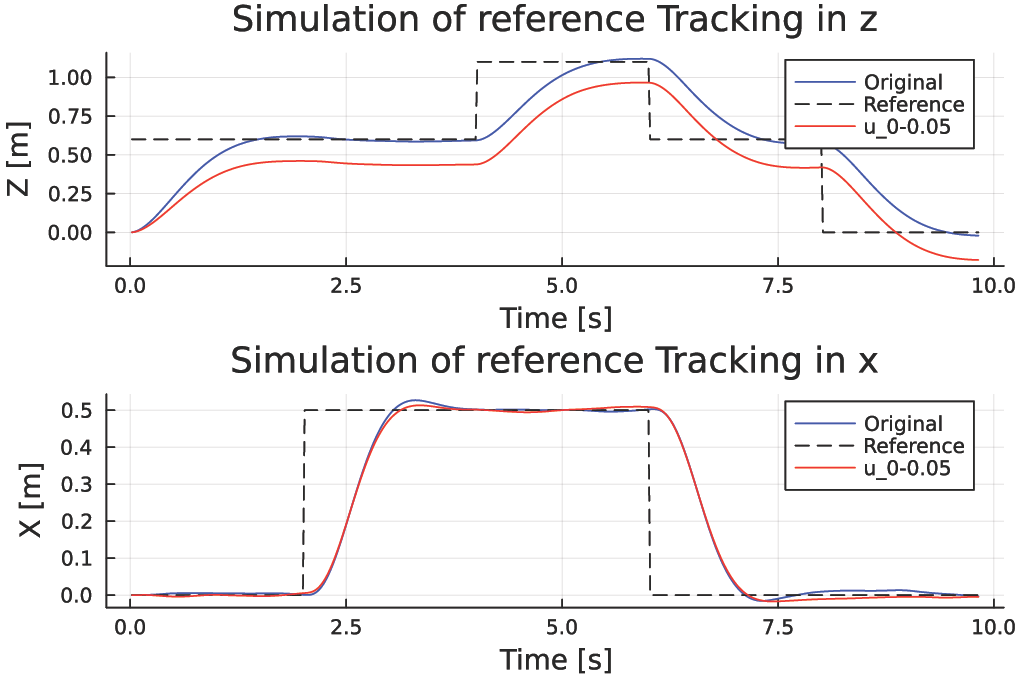}
    \caption{Simulated waypoint tracking, both with and without static offset on hover operating point $u_0$. It can be seen that unmodeled phenomena not taken into account when designing the MPC controller affect altitude control.}
    \label{fig:flightdaqptinympceightsimulation}
\end{figure}                                                                          
\end{document}